\documentclass[a4paper,12pt]{article}

\usepackage{amsmath,amsfonts,amssymb,amsthm}
\theoremstyle{plain}
\newtheorem{thm}{Theorem}
\newtheorem*{thm*}{Theorem}
\newtheorem*{Main-SAW}{Main Theorem}
\newtheorem*{TECHthm}{Technical Theorem}

\newtheorem*{lem*}{Lemma}
\newtheorem*{cor}{Corollary}

\theoremstyle{definition}

\newtheorem{Def}{Definition}

\newcommand{\reftit}{\textit}    
\newcommand{\refis}{\textbf}     

\begin{document}

\title{Brownian Bridge and Self-Avoiding Random Walk.}

\author{Yevgeniy Kovchegov\\
 \small{Email: yevgeniy@math.stanford.edu}\\
 \small{Fax: 1-650-725-4066}}

\maketitle

\begin{abstract}
We derive the Brownian bridge asymptotics for a scaled
self-avoiding walk conditioned on arriving to a far away point
 $n \vec{\mathbf{a}}$ for $\vec{\mathbf{a}} \in (\mathbb{Z},0,...,0)$,
 and outline the proof for all other $\vec{\mathbf{a}}$ in $\mathbb{Z}^d$.
\end{abstract}

\section{Introduction.}

\subsection{Self-Avoiding Walks and Brownian Bridge.}

 First we briefly introduce the notion of a self-avoiding walk based
 on the material rigorously presented in \cite{slade}, and
 the notion of the Brownian Bridge followed by the history of the
 problem solved as well as that of the tools and ideas used in producing
 the results. We will conclude the introduction with the word description
 of the results of this research paper given in light of the chronological
 perspective, with the brief outline of the proofs that follow.

 \textbf{Self-Avoiding Random Walk (SARW):}
 An $N$-step self-avoiding walk (path) $\omega$ on $\mathbb{Z}^d$,
 beginning at $0$ is a sequence of sites: $\omega(0)=0,
 \omega(1),..., \omega(N)$ with $|\omega(j+1)-\omega(j)|=1$ and
 $\omega(i) \not= \omega(j)$ for all $i \not= j$. We let $c_N$
 denote the number of $N$-step self-avoiding walks beginning at
 zero. It had been established that the limit representing the connective
 constant $\mu = \lim_{N \rightarrow \infty} {c_N}^{\frac{1}{N} }$
 exists due to a subadditivity property of $\log{c_N}$ (see
 \cite{slade}). We also let $c_N(x,y)$ to be the number of
 $N$-step self-avoiding walks $\omega$ with $\omega(0)=x$ and
 $\omega(N)=y$.  The two-point function $g_{\beta}(x,y)$ (as defined below)
 is an important tool in the theory of self-avoiding walks:
 $$g_{\beta}(x,y) \equiv \sum_{N=0}^{\infty} c_N(x,y) e^{-\beta N} =
 \sum_{\omega: x \rightarrow y} e^{-\beta |\omega|},$$
 where the second sum is taken only over all self-avoiding walks
 $\omega: x \rightarrow y$ on the lattice. For the simplicity of
 notation (due to the shift-invariance property of $c_N(x,y)$ )
 we denote $g_{\beta}(x) \equiv g_{\beta}(0,x)$. The supercritical
 $\beta > \beta_c(d)$ is the one for which the equivalent sums above are
 finite. It was shown (see \cite{slade}) that for the
 supercritical $\beta$, the "bubble diagram"
 $$B_d(\beta) \equiv \sum_{x} g_{\beta}(x)^2$$
 is finite. The significance of the bubble diagram is discussed in
 Section 1.5 of \cite{slade}.\\

 Since the radius of convergence $e^{-\beta_c(d)} =\frac{1}{\mu}$,
 it is apparent that the two-point function decays exponentially:
 \begin{eqnarray} \label{exp-decay}
 g_{\beta}(0,x) \leq C_{\beta} e^{-c_{\beta} \|x\|}
 \end{eqnarray}
 for all $\beta > \beta_c(d)$ and some corresponding $C_{\beta}, c_{\beta} >0$.

 The notion of a "mass" of a two-point function applies here as
 well. The mass $m(\beta)$ is the rate of exponential decay of $g_{\beta}(x,y)$
 in the direction of the first coordinate vector:
 $$m(\beta)= \liminf_{n \rightarrow \infty}
 \frac{-\log{g_{\beta}(0,(n,0,...,0))}}{n}. $$
 It can be shown that the $\liminf$ above can be replaced by the
 limit.

 \textbf{Brownian Bridge:} defined as a sample-continuous Gaussian
 process $B^0$ on $[0,1]$ with mean $0$ and $\bold{E} B^0_s B^0_t
 = s(1-t)$ for $0 \leq s \leq t \leq 1$. So, $B^0_0 = B^0_1 = 0$
 a.s. Also, if $B$ is a Brownian motion, then the process
 $B_t - tB_1$ ($0 \leq t \leq 1$) is a Brownian Bridge. For more
 details see \cite{bill}, \cite{dudley} and \cite{durrett}.  In a more general
 setting, we call the process
 $B^{0, \mathbf{\vec{a}}}_t \equiv B^0_t +t \vec{\mathbf{a}}$
 "a Brownian Bridge connecting points zero and $\vec{\mathbf{a}}$".\\

\subsection{History of the Problem, the Results and Strategy of the Proof.}

The main goal of this paper is to show the derivation of the weak
Brownian bridge asymptotics for a scaled self-avoiding walk
conditioned on arriving to a far away point $n \vec{\mathbf{a}}$
($\vec{\mathbf{a}} \in \mathbb{Z}^d$). The technique used in the
proofs originates from the methodology developed in the process of
establishing a precise Ornstein-Zernike decay for a variety of
spin systems and lattice field theories and the development of the
renewal theory. It turned out that the technique developed by
Ornstein and Zernike in 1914 for the case of the classical fluid
can be implemented in many classical models of statistical
mechanics (self-avoiding walks, percolation, 2D Ising model and
many other spin systems) for all noncritical temperatures. For
this, for the given two-point function, one needs to construct a
"direct correlation function" with a strictly greater rate of
decay. This approach was implemented in in the case of the
$d$-dimensional self-avoiding walks \cite{cc} giving the precise
Ornstein-Zernike behavior of the two-point function
 $g_{\beta} (0, n \mathbf{\vec{a}})$ connecting the origin to a
 point on an axis
 (the case $\mathbf{\vec{a}} = (\| \mathbf{\vec{a}} \|,0,...,0)$)
for all noncritical $\beta$. There the so called "mass gap"
condition (or separation of mass)is proved. In that case, the
two-point function with a different rate of decay is the
generating function corresponding to the self-avoiding walks with
all non-trivial (more than one) intersections with the hyper
spaces $\{x_1=c\}$ situated in between the origin and the
destination point. The work of proving the Ornstein-Zernike
behavior (the
 coefficient of order $\|x\|^{\frac{d-1}{2}}$ near the decay
 exponent of the two-point function) was completed in
 \cite{ioffe1} for any supercritical value of the parameter
 $\beta > \beta_c(d)$. There the complete precise asymptotics (\ref{OZdecay})
 of the decay was derived in any direction $\vec{\mathbf{a}}$ as the result of
 an extensive study of the geometric properties of corresponding equi-decay
 level sets, broadening the methodology of \cite{cc}.\\

 The corresponding developments in subcritical bond percolation
 model followed the above advances in the theory of self-avoiding
 walks. In \cite{acc}, \cite{ccc} and \cite{ioffe} some similar
 equi-decay level sets are studied, and corresponding Ornstein-Zernike
 asymptotics is produced. This technique was used in
 \cite{kovcheg} together with the technical result of section \ref{intro:2} to
 produce a Brownian bridge asymptotics of a scaled percolation
 cluster conditioned on reaching a far away point, and also proving
 the shrinking of such clusters. In this paper, we follow up on
 the result of \cite{kovcheg}. We prove the weak convergence of a
 scaled interpolation "skeleton" going through the regeneration points
 (see definition \ref{reg:pt})
 of a self-avoiding walk, and terminating at a far away point $n\vec{\mathbf{a}}$
 to Time$\times$($d-1$)-dimensional Brownian bridge as $n \rightarrow \infty$.
 Later, the shrinking of the self-avoiding walk to the above interpolation
 skeleton is proved (see section \ref{shrink:SAW}). We prove the result for
  $\mathbf{\vec{a}} = (\| \vec{\mathbf{a}} \|,0,...,0)$ given an
 appropriate measure on such self-avoiding walks (see (\ref{meas:walk})).
 We outline the proof of the result for all other $\vec{\mathbf{a}}$ in
 $\mathbb{Z}^d$.

\subsection{Asymptotic Convergence to Brownian Bridge.}\label{intro:2}

 The following technical result was proved in \cite{kovcheg}.
 Let $X_1, X_2,...$ be i.i.d. random variables on $\mathbb{Z}^d$
  with the span of the lattice distribution equal to one
(see \cite{durrett}, section 2.5), and let there be a
 $\bar{\lambda} > 0$ such that the moment-generating function
 $$\bold{E}(e^{\theta \cdot X_1}) <\infty$$
for all $\theta \in B_{\bar{\lambda}}$.
\\

 Now, for a given vector $\mathbf{\vec{a}} \in \mathbb{Z}^d$,
 let $X_1+...+X_i =[t_i, Y_i]_f \in \mathbb{Z}^d$ when
 written in the new orthonormal basis such that
 $\mathbf{\vec{a}} = [\| \mathbf{\vec{a}} \|, 0]_f$ (in the new basis
 $[\cdot , \cdot ]_f \in \mathbb{R} \times \mathbb{R}^{d-1}$).
 Also let $P[\mathbf{\vec{a}} \cdot X_i]>0]=1$.
 We define the process $[t, Y_{n,k}^*(t)]_f$ to be
 the interpolation of $0$ and
 $[\frac{1}{n \| \mathbf{\vec{a}} \|}t_i, \frac{1}{\sqrt{n}}Y_i]_f^{i=0,1,...,k}$,
 in Section 2.2 we will show that

\begin{TECHthm}
The process
 \begin{eqnarray} \label{tech}
 \{ Y^*_{n,k}\mbox{ for some } k \mbox{ such that }[t_k, Y_k]_f =n \mathbf{\vec{a}} \}
 \end{eqnarray}  conditioned on
the existence of such $k$ converges  weakly to the
 Brownian Bridge (of variance  that depends only on the law of $X_1$).
\end{TECHthm}

\section{The Main Result in SARW.} \label{SARW}

 In this section we work only with supercritical SARW ($\beta > \beta_c(d)$).

\subsection{Preliminaries.} \label{SAW:prelim}

Here we briefly go over the definitions that one can find in
\cite{slade}. We start with the decay rate
$\tau_{\beta}(\vec{x})$:

$$\tau_{\beta}(\vec{x}) \equiv -\lim_{n \rightarrow \infty}
   \frac{1}{n} \log{g_{\beta}([n\vec{x}])},$$
where the limit is always defined since
 $$\frac{g_{\beta}(\vec{x}+\vec{y})}{B_d(\beta)} \geq
  \frac{g_{\beta}(\vec{x})}{B_d(\beta)} \frac{g_{\beta}(\vec{y})}{B_d(\beta)}.$$
 Now, $\tau_{\beta}(\vec{x})$ is the support function of the compact convex set
$$ \bold{K}^{\beta} \equiv \bigcap_{\vec{n} \in \mathbb{S}^{d-1}}
  \lbrace \vec{r} \in \mathbb{R}^d \mbox{ : }
  \vec{r} \cdot \vec{n} \leq \tau_{\beta}(\vec{n}) \rbrace ,$$
with non-empty interior int\{$\bold{K}^{\beta}$\} containing point
zero.
\\
Let $\omega(j)=(\omega_1(j),...,\omega_d(j))$ be a self-avoiding
path defined for $j \in [a,b] \bigcap \mathbb{N}$, $a \leq b \in
\mathbb{Z}^+$.

\begin{Def}
We call $\omega$ \textbf{a bridge} if
 $$\omega_1(a)< \omega_1(j) \leq \omega_1(b)$$
for all $a<j \leq b$. If $x=\omega(a)$ is the initial point and
$y=\omega(b)$ is the final point, we write
 $\omega \mbox{ : } x -^b\rightarrow y$.
\end{Def}

 For $\vec{x} \in \mathbb{Z}^d$, we define the \textit{cylindrical}
 two-point function
 $$h(\vec{x}) \equiv \sum_{\omega : 0 -^b\rightarrow \vec{x}}
 { e^{-\beta |\omega|} } ,$$
where  $h(\vec{x}) = \delta_0 (\vec{x})$ for all
 $\vec{x} \in \{0\} \times \mathbb{Z}^{d-1}$ .

\begin{Def}
We say that $k \in \mathbb{N}$ ($\omega_1(a)<k<\omega_1(b)$) is
\textbf{a break point} of $\omega$ if there exists $r \in [a,b]$
such that $\omega_1(j) \leq k$ whenever $j \leq r$ and
$\omega_1(j)>k$ whenever $j>r$.
\end{Def}

\begin{Def}
A bridge $\omega \mbox{ : } x -^b\rightarrow y$ (where, as before,
$x=\omega(a)$ and $y=\omega(b)$) is called \textbf{irreducible} if
it has no break points. In that case we write $\omega \mbox{ : } x
-^{ib}\rightarrow y$.
\end{Def}

Now, for $\vec{x} \in \mathbb{Z}^d$, we define the
\textit{irreducible}
 two-point function
 $$f(\vec{x}) \equiv \sum_{\omega : 0 -^{ib}\rightarrow \vec{x}}
 { e^{-\beta |\omega|} } ,$$
with  $f(\vec{x}) = \delta_0 (\vec{x})$ for all
 $\vec{x} \in \{0\} \times \mathbb{Z}^{d-1}$ .

\subsection{SARW and Regeneration Structures.} \label{SAW:reg}

It turned out that if counting the bridges between the origin and
a point $\vec{k}=(k_x, k_y) \in \mathbb{N} \times
\mathbb{Z}^{d-1}$, that $f$ and $h$ satisfy the recurrence
equation (see \cite{slade}):
\begin{eqnarray} \label{OZ-SAW}
 h(\vec{k})= \sum_{i=1}^{k_x} \sum_{l \in \mathbb{Z}^{d-1}} f(i,l)h(k_x -i, k_y -l),
\end{eqnarray}
 which together with
$h(0,\tilde{k})=\delta_0(\tilde{k})$ (for $\tilde{k} \in
\mathbb{Z}^{d-1}$) are called the Ornstein-Zernike equations.\\
 Now, for any $\tilde{r} \in \mathbb{Z}^{d-1}$, we define
\begin{eqnarray} \label{FH}
   H_n(\tilde{r}) \equiv \sum_{\tilde{k} \in \mathbb{Z}^{d-1}}
   e^{\tilde{r} \cdot \tilde{k}} h(n,\tilde{k}) \mbox{   and   }
   F_n(\tilde{r}) \equiv \sum_{\tilde{k} \in \mathbb{Z}^{d-1}}
   e^{\tilde{r} \cdot \tilde{k}} f(n,\tilde{k}),
\end{eqnarray}
as well as the corresponding mass
 $$m_H(\tilde{r}) \equiv \lim_{n \rightarrow +\infty} \frac{1}{n}
   \log{H_n(\tilde{r})}  \mbox{   and   }
   m_F(\tilde{r}) \equiv \lim_{n \rightarrow +\infty} \frac{1}{n}
   \log{F_n(\tilde{r})}.$$

The Ornstein-Zernike asymptotics has been proved for the
cylindrical two-point function $h(\cdot)$ (see \cite{cc} and
\cite{ioffe1}), using the "mass gap" condition, e.g. existence of
a point $\tilde{r}_o \in \mathbb{Z}^{d-1}$, inside a neighborhood
of points with finite mass $m_H$, such that $m_H(\tilde{r}_o) >
m_F(\tilde{r}_o)$. It was also shown (see \cite{ioffe1}, Section
2) that the mass gap condition with the renewal theorem
(\cite{slade}, Appendix B) imply that
$\exp\{-nm_H(\tilde{r}_o)\}F_n(\tilde{r}_o)$ is a probability
distribution (where $\tilde{r}_o$ is as above):
\begin{eqnarray} \label{dist}
 \sum_{n \in \mathbb{N}}\exp\{-nm_H(\tilde{r}_o)\}F_n(\tilde{r}_o)=1.
\end{eqnarray}

As it was mentioned in the introduction, the mass gap condition
was crucial in obtaining the Ornstein-Zernike decay (see
\cite{ioffe1}):
\begin{thm} \label{OZdecay}
For all $d \geq 2$ and $\beta > \beta_c(d)$,
$$g_{\beta}(x)= \psi_{\beta}(\frac{x}{\|x\|})
  \frac{e^{-\tau_{\beta}(x)}}{\|x\|^{\frac{d-1}{2}}} (1+o(1))$$
uniformly in $\|x\|$, where $\psi_{\beta}(\cdot)$ is analytic on
the unit circle.
\end{thm}

\subsection{Measure $Q_{r_0}(x)$.}
We notice that substituting the sum $F_n(\tilde{r}_o)$, as defined
in (\ref{FH}), into (\ref{dist}) we obtain (after some simple
manipulations) an enhanced version of (\ref{dist}):
$$\sum_{\vec{x} \in \mathbb{N} \times \mathbb{Z}^{d-1}} f(\vec{x})
  e^{\vec{x} \cdot (-m_H(\tilde{r}_o),\tilde{r}_o)} =1,$$
where $\vec{r}_o \equiv (-m_H(\tilde{r}_o),\tilde{r}_o) \in
\partial \bold{K}^{\beta}$ as it was shown in \cite{ioffe1}, Section 3.
\\

Now, let for $\vec{x} \in \mathbb{N} \times \mathbb{Z}^{d-1}$,
$$Q_{r_0}(\vec{x}) \equiv f(\vec{x}) e^{\vec{x} \cdot \vec{r}_0}.$$
Due to the equation above, $Q_{r_0}(\cdot)$ is a probability
measure on $\mathbb{N} \times \mathbb{Z}^{d-1}$. It is similar to
the regeneration measure, defined for the subcritical bond
percolation model in Section 4 of \cite{ioffe}, and later used in
\cite{kovcheg} for derivation of Brownian Bridge assymptotics for
that model.\\

 The mass gap condition implies the exponential decay of $Q_{r_0}(\vec{x})$.

\subsection{The Result For $\mathbf{\vec{a}}=(1,0,...,0)$.}

 We fix $\bold{\vec{a}} \in \mathbb{Z}^d$. We let for a supercritical
 constant $\beta$ and all $n \in \mathbb{N}$,
 $P_n(\cdot)$ to be a
 law on a set of self-avoiding random paths $\omega$, conditioned on $\omega$
 being a bridge between $0$ and $n \mathbf{\vec{a}}$
 ($\omega : 0 -^b\rightarrow n\mathbf{\vec{a}}$). More precisely,
 we define $P_n$ as
 \begin{eqnarray} \label{meas:walk}
 P_n(\omega : 0 -^b\rightarrow n\mathbf{\vec{a}}) \equiv
   \frac{\exp{(-\beta |\omega|)} }
    { \sum_{\tilde{\omega} : 0 -^b\rightarrow n\mathbf{\vec{a}}}
       \exp{(-\beta |\tilde{\omega}|)} }
   = \frac{\exp{(-\beta |\omega|)} }{h(n\mathbf{\vec{a}})}.
 \end{eqnarray}
For now, we let $\mathbf{\vec{a}}=[1,0] \equiv (1,0,...,0)$ and
 $\vec{r}_o = (\mathbb{Z}^+,0,...,0) \bigcap \partial \bold{K}^{\beta}$.
Here, we define \textit{the regeneration points} in a way, similar
to that, used to define the regeneration points for the case of
Bernoulli bond percolation model:
\begin{Def} \label{reg:pt}
Given a self-avoiding path $\omega$, and a break point $b$. We say
that $\omega(T_b)$ is the corresponding \textbf{regeneration
point} if $T_b=\max\{t \mbox{ : } \omega_1(t)=b \}$.
\end{Def}
In a very important step, we notice that given the
Ornstein-Zernike equations (\ref{OZ-SAW}) above and the definition
of probability distribution $P_n(\cdot)$, we can explicitly write
(in terms of measure $Q_{r_o}$) the probability of the walk
passing through the particular regeneration points $s_1 \equiv
x_1$, $(s_2 \equiv  x_1+x_2)$,..., $(s_{k-1} \equiv
x_1+...+x_{k-1})$, where are all $x_i \in \mathbb{Z}^+ \times
\mathbb{Z}^{d-1}$:
\begin{eqnarray} \label{imp}
P_n[s_1,...,s_{k-1} \mbox{  are reg. pts.  }] & = &
 \frac{1}{h(n\mathbf{\vec{a}})}
 \left(\sum_{\omega : 0 -^{ib}\rightarrow s_1}{e^{-\beta|\omega|}}\right)
 ...
 \left(\sum_{\omega :s_{k-1} -^{ib}\rightarrow s_k}{e^{-\beta|\omega|}}\right) \nonumber \\
 & = &
 \frac{1}{h(n\mathbf{\vec{a}})}f(x_1)...f(x_k) \nonumber \\
 & = &
 \frac{Q_{r_o}(x_1)...Q_{r_o}(x_k)}
      {\sum_{\kappa} \sum_{\varsigma_1+...+\varsigma_{\kappa}=n\mathbf{\vec{a}}}
       Q_{r_o}(\varsigma_1)...Q_{r_o}(\varsigma_{\kappa})},
\end{eqnarray}
where $s_0 \equiv 0$ and $s_k \equiv x_1+...+x_{k-1}+x_k =
n\mathbf{\vec{a}}$.
\\

 We recall that the moment generating function (the Laplace
 transform) under the measure $Q_{r_o}(\cdot)$ is finite in a
 neighborhood of zero:
 $$\bold{E}_{r_0} (e^{\theta \cdot X_1}) < \infty$$
 for all $\theta \in B_{\bar{\lambda}}(0)$.
 We use the brackets
 $[\cdot,\cdot] \in \mathbb{R} \times \mathbb{R}^{d-1}$ to denote
 the coordinates of $\mathbb{R}^d$ vectors for the simplicity of
 notation.
 Obviously $\mathbf{\vec{a}}=[1, 0]$.
 We want to prove that the process corresponding to the last $d-1$
coordinates in the new basis of the scaled
 ($\frac{1}{n}$ times along $\bold{\vec{a}}$
 and $\frac{1}{\sqrt{n}}$ times in the orthogonal $d-1$
 dimensions) interpolation of regeneration points
  of the self-avoiding path $\omega$  conditioned on
 $\omega : 0 -^b\rightarrow n\mathbf{\vec{a}}$
 converges weakly to the Brownian Bridge $B^o(t)$ (with variance that depends
 only on measure $Q_{r_0}$) where $t$ represents the scaled down first coordinate.\\

Let $X_1, X_2,...$ be i.i.d. random variables distributed
according to $Q_{r_0}$ law.
 We interpolate $0,X_1,(X_1+X_2),...,(X_1+...+X_k)$ and scale by
$\frac{1}{n} \times \frac{1}{\sqrt{n}}$ along $<\mathbf{\vec{a}}>
\times <\mathbf{\vec{a}}>^{\bot}$ to get the process $[t,
Y_{n,k}(t)]$. The technical theorem (\ref{tech}) implies the
following

\begin{thm}
The process
$$\{ Y_{n,k}\mbox{ for some } k \mbox{ such that }
    X_1+...+X_k= n \bold{\vec{a}} \}$$
 conditioned on  the existence of such $k$ converges weakly to
 the Brownian Bridge (with variance that depends only on measure $Q_{r_0}$).
\end{thm}

Now, let for $y_1,...,y_k \in \mathbb{Z}^d$ with  positive
increasing first coordinates,
 $\gamma (y_1,...,y_k)$  be the last $(d-1)$ coordinates
  in the new basis
 of the scaled ($\frac{1}{n} \times \frac{1}{\sqrt{n}}$)
interpolation of points $0,y_1,...,y_k$ (where the first
coordinate is time). Notice that $\gamma (y_1,...,y_k) \in
C_o[0,1]^{d-1}$
 as a function of scaled first coordinate
 whenever $y_k= n \bold{\vec{a}}$.
\\
 By the important observation (\ref{imp}) we've made before,
 for any function $F(\cdot )$ on $C[0,1]^{d-1}$,\\
\\
$\sum_k \sum_{x_1+...+x_k= n \bold{\vec{a}}}
  F(\gamma (x_1, x_1+x_2, ... , \sum_{i=1}^k x_i))$
$$ \times
 P_n[0 \leftarrow^{h_{r_0}} \rightarrow x
  \mbox{ ;  regeneration points:  }
  x_1, x_1+x_2, ... , \sum_{i=1}^k x_i =x]$$

$$=\sum_k \sum_{x_1+...+x_k= n \bold{\vec{a}}}
 F(\gamma (x_1, x_1+x_2, ... , \sum_{i=1}^k x_i))
 f(x_1)...f(x_k)$$

$$= e^{-r_0 \cdot n \bold{\vec{a}}}
 \sum_k \sum_{x_1+...+x_k= n \bold{\vec{a}}}
 F(\gamma (x_1, x_1+x_2, ... , \sum_{i=1}^k x_i))
 Q_{r_0}(x_1)...Q_{r_0}(x_k).$$
\\
Therefore, for any $A \subset C[0,1]^{d-1}$\\
\\
$P_p[ \gamma (\mbox{regeneration points of }\omega) \in A \mbox{ }
| \mbox{  }\omega : 0 -^b\rightarrow n\mathbf{\vec{a}} ]$

$$
=\frac{  \sum_k \sum_{x_1+...+x_k= n \bold{\vec{a}}}
         I_A (\gamma (x_1, x_1+x_2, ... , \sum_{i=1}^k x_i)) f(x_1)...f(x_k) }
        { \sum_k \sum_{x_1+...+x_k= n \bold{\vec{a}}} f(x_1)...f(x_k) }$$

$$
=\frac{\sum_k \sum_{x_1+...+x_k= n \bold{\vec{a}}}
         I_A (\gamma (x_1, x_1+x_2, ... , \sum_{i=1}^k x_i))
         Q_{r_0}(x_1)...Q_{r_0}(x_k)}
        {\sum_k \sum_{x_1+...+x_k= n \bold{\vec{a}}}
         Q_{r_0}(x_1)...Q_{r_0}(x_k)}$$
$$= P[Y_{n,k} \in A \mbox{ for the } k \mbox{ such that }
    X_1+...+X_k= n \bold{\vec{a}} \mbox{  } | \mbox{  }
  \exists k \mbox{ such that }
    X_1+...+X_k= n \bold{\vec{a}}] .$$

Hence, we have proved the following

\begin{cor}
 The process corresponding to the last $d-1$ coordinates of the scaled
 $({\frac{1}{n } \times \frac{1}{\sqrt{n}} })$
 interpolation of regeneration points of the self-avoiding path
 $\omega$ (with the scaled first coordinate denoting the time interval) conditioned on
 $\omega : 0 -^b\rightarrow n\mathbf{\vec{a}}$
 converges weakly to the Brownian Bridge (with variance that depends
 only on measure $Q_{r_0}$).
\end{cor}

\subsection{Shrinking of the Self-Avoiding Walks.} \label{shrink:SAW}
 Here we again let $\mathbf{\vec{a}}=[1,0] \equiv (1,0,...,0)$ and
 $\vec{r}_o = [\mathbb{Z}^+,0] \bigcap \partial \bold{K}^{\beta}$.
 In the way of proving that the scaled walk
 $\omega : 0-^b\rightarrow n\mathbf{\vec{a}}$ shrinks, we shell need
 to show that the consequent
 regeneration points are situated relatively close to each other:
\begin{lem*}
$$P_p[\max_{i} |x_i - x_{i-1}|>n^{1/3},
 \mbox{  } x_i \mbox{- reg. points  }
  | \mbox{  } 0-^b\rightarrow n\mathbf{\vec{a}} ]
 <\frac{1}{n}$$
for $n$ large enough.
\end{lem*}

\begin{proof}
 Since $\vec{r}_o =[\| \vec{r}_o \|, 0] \in \partial \bold{K}^{\beta}$
 and therefore
 $$\vec{r}_o \cdot [v_1,v_2] = \vec{r}_o \cdot [v_1, 0] \leq
 \tau_{\beta}([v_1,0]) \leq \tau_{\beta}([v_1, v_2])$$
 for all $[v_1,v_2] \in Z^d$,
 $$\vec{r}_o \cdot [1,0] = \tau_{\beta}([1,0]) \mbox{    and     }
 \| \vec{r}_o \|^2 = \tau_{\beta}(\vec{r}_o).$$
 Hence, by the pseudo-linearity of $\tau_{\beta}(\cdot)$,
 $$\nabla \tau_{\beta}(\bold{\vec{a}})=[\tau_{\beta}(1),0]
   =[\| \vec{r}_o \|, 0] = \vec{r}_o.$$

Now, by the convexity of $\tau_{\beta}(\cdot)$,
$$ \frac{ \tau_{\beta}(\bold{\vec{a}}) - \tau_{\beta}(\bold{\vec{a}} -\frac{\vec{x}}{n})}
 {(\frac{\|\vec{x}\|}{n})}
 \leq \frac{\vec{x}}{\|\vec{x}\|} \cdot \nabla \tau_{\beta}(\bold{\vec{a}}) $$
for $\vec{x} \in \mathbb{Z}^d$ ($\vec{x} \not= 0$), and therefore
$$\tau_{\beta}(n \bold{\vec{a}})- \tau_{\beta}(n \bold{\vec{a}}-\vec{x})
 = \|\vec{x}\| \frac{ \tau_{\beta}(\bold{\vec{a}}) - \tau_{\beta}(\bold{\vec{a}}
-\frac{\vec{x}}{n})}{(\frac{\|\vec{x}\|}{n})}
 \leq \vec{x} \cdot \nabla \tau_{\beta}(\bold{\vec{a}})
 = \mathbf{\vec{r}} \cdot \vec{x}.$$
Thus, since $Q_{r_0}(\vec{x})$ decays exponentially and therefore
$$f(\vec{x}) e^{\tau_{\beta}(n \bold{\vec{a}})- \tau_{\beta}(n \bold{\vec{a}}-\vec{x})}
  \leq Q_{r_0}(\vec{x})$$
and also decays exponentially. Hence by Ornstein-Zernike result
(Theorem \ref{OZdecay}),

$$P_p[  n^{1/3} <
  |\vec{x}| , \mbox{ } \vec{x} \mbox{-first reg. point  }
  | 0-^b\rightarrow n\mathbf{\vec{a}} ]
 =\sum_{ n^{1/3} < |\vec{x}| }
 f(\vec{x}) \frac{h(n \bold{\vec{a}}-\vec{x})} {h(n \bold{\vec{a}})}
 < \frac{1}{n^2} $$
for $n$ large enough.
 So, since the number of the regeneration points is no greater
 than $n$,

 $$P_p[\max_{i} |x_i - x_{i-1}|>n^{1/3},
 \mbox{  } x_i \mbox{- reg. points  }
  | \mbox{  } 0-^b\rightarrow n\mathbf{\vec{a}} ]
 <\frac{1}{n}$$
for $n$ large enough.

\end{proof}

 Now,  it is really easy to check that there is a
constant $\lambda_f >0$ such that
 $$f(\vec{x}) > e^{-\lambda_f \|\vec{x}\|}$$
for all $\vec{x}$ such that $f(\vec{x}) \not= 0$ (here we only
need to connect points zero and $\vec{x}$ with an "S"-shaped
irreducible bridge). Hence, due to the exponential decay
(\ref{exp-decay}) of the two point function $g_{\beta}(x,y)$, for
a given $\epsilon >0$,
$$P_p[\mbox{ the walk } \{ \omega(i) \}_{i=0,...,|\omega(\vec{x})|}
 \not\subset [\mathbb{R}, B_{\epsilon \sqrt{n}}^{d-1}(0)]
 \mbox{  } | \mbox{  } 0-^{ib}\rightarrow \vec{x}]
  < C_{\beta} e^{ \lambda_f \|\vec{x}\| -c_{\beta} \epsilon \sqrt{n}}, $$
and therefore, summing over the regeneration points, we get

$$P_p[\mbox{ the scaled walk } \{ \omega(i) \}_{i=0,...,|\omega(\vec{x})|}
 \not\subset \epsilon \mbox{-neighbd. of }
 [0,1] \times \gamma (\mbox{ reg. points })
 \mbox{  } | \mbox{  } 0-^b\rightarrow n\mathbf{\vec{a}} ]$$

$$ < \frac{1}{n} +
 n C_{\beta} e^{ \lambda_f \|\vec{x}\| -c_{\beta} \epsilon \sqrt{n}} $$
 for $n$ large enough due to the lemma above.\\

We can now state the main result for $\mathbf{\vec{a}}=[1,0]$:
\begin{Main-SAW}
 The process corresponding to the last $d-1$ coordinates of the scaled
 $({\frac{1}{n } \times \frac{1}{\sqrt{n}} })$
 interpolation of regeneration points of the self-avoiding path
 $\omega$ (with the scaled first coordinate denoting the time interval) conditioned on
 $\omega : 0 -^b\rightarrow n\mathbf{\vec{a}}$
 converges weakly to the Brownian Bridge (with variance that depends
 only on measure $Q_{r_0}$).
\\
 Also for a given $\epsilon >0$
$$P_p[\mbox{ the scaled walk } \{ \omega(i) \}_{i=0,...,|\omega(\vec{x})|}
 \not\subset \epsilon \mbox{-neighbd. of }
 [0,1] \times \gamma (\mbox{ reg. points })
 \mbox{  } | \mbox{  } 0-^b\rightarrow n\mathbf{\vec{a}} ]
 \rightarrow 0$$
as $n \rightarrow \infty$.
\end{Main-SAW}

\subsection{General Case.}

 Now, we turn our attention to all $\vec{\mathbf{a}} \in \mathbb{Z}^d$ not
 on the axis.
 It turned out that the main theorem of section \ref{shrink:SAW}
 holds for all $\vec{\mathbf{a}}$ in $\mathbb{Z}^d$. In a more
 direct approach used in the corresponding developments in percolation
 (see Section 4 of \cite{ioffe}) and finite range Ising models (see \cite{ioffe2}),
 we can replicate the same recurrence structures, as those in section
 \ref{SAW:prelim}, in a given direction (say $\vec{\mathbf{a}}$), yielding
 the same renewal relations (as in section \ref{SAW:reg}). The technique is
 simpler than that used in percolation and finite range Izing models. We
 choose a direction vector $\vec{\mathbf{r}} \in \partial \mathbf{K}^{\beta}$,
 we define the corresponding notions of "a bridge" in the
 direction $\vec{\mathbf{r}}$ and the $\vec{\mathbf{r}}$-regeneration points:

 \begin{Def}
 We call $\omega$ \textbf{an $\vec{\mathbf{r}}$-bridge} if
 $$\omega(a) \cdot \vec{\mathbf{r}} < \omega(j) \cdot \vec{\mathbf{r}}
   \leq \omega(b) \cdot \vec{\mathbf{r}}$$
 for all $a<j \leq b$. If $x=\omega(a)$ is the initial point and
 $y=\omega(b)$ is the final point, we write
 $\omega \mbox{ : } x -^{b(\vec{\mathbf{r}})}\rightarrow y$.
 \end{Def}

 Similarly, we define the \textit{cylindrical}
 two-point function
 $$h_{\vec{\mathbf{r}}}(\vec{x}) \equiv
 \sum_{\omega : 0 -^{b(\vec{\mathbf{r}})}\rightarrow \vec{x}}{ e^{-\beta |\omega|} },$$
 where  $h_{\vec{\mathbf{r}}}(\vec{x}) = \delta_0 (\vec{x})$ for all
 $\vec{x} \in <\vec{\mathbf{r}}>^{\perp}$.

 \begin{Def}
 We say that $\omega(k) \in \mathbb{Z}^d$ ($a<k<b$) is
 \textbf{an $\vec{\mathbf{r}}$-regeneration point} of $\omega$ if there exists
 $N \in [a,b]$ such that
 $\omega(j) \cdot \vec{\mathbf{r}} \leq \omega(k) \cdot \vec{\mathbf{r}}$
 whenever $j \leq N$ and
 $\omega(j) \cdot \vec{\mathbf{r}} > \omega(k) \cdot \vec{\mathbf{r}}$
 whenever $j>N$.
 \end{Def}

 \begin{Def}
 An $\vec{\mathbf{r}}$-bridge $\omega \mbox{ : } x -^b\rightarrow y$
 (where, as before, $x=\omega(a)$ and $y=\omega(b)$) is called
 $\omega(k) \cdot \vec{\mathbf{r}}$-\textbf{irreducible} if it has
 no $\vec{\mathbf{r}}$-regeneration points. In that case we write
 $\omega \mbox{ : } x-^{ib(\vec{\mathbf{r}})}\rightarrow y$.
 \end{Def}

 We again redefine the corresponding \textit{irreducible}
 two-point function
 $$f_{\vec{\mathbf{r}}}(\vec{x}) \equiv
 \sum_{\omega : 0 -^{ib(\vec{\mathbf{r}})}\rightarrow \vec{x}}{ e^{-\beta |\omega|} },$$
 with  $f_{\vec{\mathbf{r}}}(\vec{x}) = \delta_0(\vec{x})$ for all
 $\vec{x} \in <\vec{\mathbf{r}}>^{\perp}$.\\

 The generalized Ornstein-Zernike recurrence equations also hold
 here: by counting the $\vec{\mathbf{r}}$-bridges between the origin
 and a lattice point
 $\vec{k} \in \mathbb{<\vec{\mathbf{r}}>} \times <\vec{\mathbf{r}}>^{\perp}$,
 we establish
 \begin{eqnarray}
 h_{\vec{\mathbf{r}}}(\vec{k})=
 \sum_{0<\vec{m} \cdot \vec{\mathbf{r}} \leq \vec{k} \cdot \vec{\mathbf{r}}}
 f_{\vec{\mathbf{r}}}(\vec{m})h_{\vec{\mathbf{r}}}(\vec{k}-\vec{m}),
 \end{eqnarray}
 where, in the sum, all $\vec{m} \in \mathbb{Z}^d$.\\

 As in \cite{ioffe}, we can replicate all the regeneration
 structures, and in particular show the existence of a positive
 $\bar{\lambda}$ such that
 $$Q^{\vec{\mathbf{r}}}_{r_0}(\vec{x}) \equiv
 f_{\vec{\mathbf{r}}}(\vec{x}) e^{\vec{x} \cdot \vec{r}_0}$$
 is a probability measure whenever $\vec{r}_0 \in
 B_{\bar{\lambda}}(\vec{\mathbf{r}}) \bigcap \partial
 \mathbf{K}^{\beta}$.
 Taking an appropriate $\vec{\mathbf{r}}$ (say
 $\vec{\mathbf{r}}=\triangledown \tau_{\beta}(\vec{\mathbf{a}})$), we can show,
 as it was done in \cite{kovcheg} for percolation clusters in subcritical phase,
 the skeleton convergence and shrinking of the scaled self-avoiding walks,
 conditioned on arriving to $n \vec{\mathbf{a}}$. Whence the main
 theorem of section \ref{shrink:SAW} would hold if we scale the
 walks by $\frac{1}{n \| \vec{\mathbf{a}} \|}$ along
 $<\vec{\mathbf{a}}>$ and by $\frac{1}{\sqrt{n}}$ in all
 orthogonal directions (along $<\vec{\mathbf{a}}>^{\perp}$).

\section*{Acknowledgements}

 The author wishes to thank A.Dembo and D.Ioffe for providing
 him with valuable and insightful comments and suggestions concerning the
 matter of this research.

\bibliographystyle{amsplain}

\end{document}